\documentclass[10pt]{article}
\usepackage{amssymb, amsmath, url, graphicx, setspace, geometry}
\usepackage{calrsfs}
\usepackage{wasysym}

\def\3{\subset }
\def\4{\subseteq }
\def\<{\left<}
\def\>{\right>}

\def\bit{\begin{itemize}}
\def\eit{\end{itemize}}
\def\3{\subset }
\def\4{\subseteq }

\def\0{\leqno}

\def\barr{\begin{array}}
\def\earr{\end{array}}

\def\Z{{\rlap{$\kern2pt{\rm Z}$}{\rm Z}\,}}
\def\bld#1#2{{\buildrel{#1}\over{#2}}}
\def\st#1#2{{\mathrel{\mathop{#2}\limits_{#1}}{}\!}}
\def\stb#1#2#3{{\st{{#1}}{\bld{{#2}}{#3}}{}\!}}
\def\xmare#1#2{\stb{#1}{#2}{\mbox{\Huge$\times$}}}

\title{\bf A connection between the number of subgroups and the order of a finite group}
\author{Mihai-Silviu Lazorec}
\date{January 18, 2019}

\begin{document}

\maketitle

\begin{abstract}
For a finite group $G$, we associate the quantity $\beta(G)=\frac{|L(G)|}{|G|}$, where $L(G)$ is the subgroup lattice of $G$. Different properties and problems related to this ratio are studied throughout the paper. We determine the second minimum value of $\beta$ on the class of $p$-groups of order $p^n$, where $n\geq 3$ is an integer. We show that the set containing the quantities $\beta(G)$, where $G$ is a finite (abelian) group, is dense in $[0,\infty).$ Finally, we consider $\beta$ to be a function on $L(G)$ and we mark some of its properties, the main result being the classification of finite abelian $p$-groups $G$ satisfying $\beta(H)\leq 1, \ \forall \ H\in L(G).$  
\end{abstract}



\section{Introduction}

One of the main characteristics of a finite group $G$ is its subgroup lattice. We denote it by $L(G)$. The connections between $G$ and $L(G)$ constitute a fruitful research topic. In this regard, some interesting problems, that were studied in the last decades, are mentioned in the Preface of the monograph \cite{21}, which is one of the well-known references on the subject. Also, the determination of the quantity $|L(G)|$ , where $G$ belongs to a remarkable class of finite groups, is a problem that is still frequently studied.    As we will remark, a lot of results on this matter were obtained especially for finite $p$-groups. In this paper, we also focus on the quantity $|L(G)|$, but we relate it to $|G|$. Therefore, for a finite group $G$, we study the quantity
$$\beta(G)=\frac{|L(G)|}{|G|}.$$

Our paper is organized as follows. As a starting point, in Section 2, we indicate some of the most relevant results concerning the counting of subgroups of a finite (abelian) $p$-group and we point out some basic properties of $\beta$. In Section 3, we determine the second minimum value of $\beta$ on the class of finite $p$-groups of order $p^n$, where $n\geq 3$ is a positive integer. As a consequence, we classify all finite $p$-groups $G$ or order $p^n$, with $n\geq 3$, that satisfy $\beta(G)\leq q_{p,n}$, where $q_{p,n}$ is a quantity that depends on $p$ and $n$.
In Section 4, we prove that the sets $\lbrace \beta(G) \ | \ G\in\mathcal{A}\rbrace$ and $\lbrace \beta(G) \ | \ G\in\mathcal{F}\rbrace$ are dense in $[0,\infty)$, where $\mathcal{A}$ is the class of finite abelian groups, while $\mathcal{F}$ is the class of all finite groups. In Section 5, instead of studying $\beta$ on a class of finite groups, we choose to view it as a function on the subgroup lattice of a finite group $G$. As an application, we classify all finite abelian $p$-groups $G$ satisfying $\beta(H)\leq 1$, for any subgroup $H$ of $G$. Finally, some further research directions are indicated in Section 6.          
 
Most of our notation is standard and will usually not be repeated here. We only mention that we denote by $L_1(G)$ the poset of cyclic subgroups of a finite group. Also, we recall that the generalized quaternion group $Q_{2^n}$, where $n\geq 3$ is an integer, and the modular $p$-groups $M(p^n)$, where $p$ is a prime and $n$ is an integer such that $n\geq 3$ if $p\geq 3$, or $n\geq 4$ if $p=2$, have the following structures:
$$Q_{2^n}=\langle x,y \ | \ x^{2^{n-1}}=y^4=1, yxy^{-1}=x^{2^{n-1}-1}\rangle \text{ and } M(p^n)=\langle x, y \ | \ x^{p^{n-1}}=y^p=1, x^y=x^{1+p^{n-2}}\rangle.$$
Elementary notions and results on groups can be found in \cite{22}. For subgroup lattice concepts we refer the reader to \cite{21}.   

\section{Preliminary results}

\subsection{Counting subgroups of finite (abelian) $p$-groups}

Let $G$ be a $p$-group of order $p^n$, where $n$ is a positive integer. For each integer $k$ such that $0\leq k\leq n$, we denote the number of subgroups of order $p^k$ of $G$ by $s_k(G)$. Two remarkable results involving these quantities are the following ones:\\

\textbf{Theorem 2.1.1.} (see Section 4 of \cite{12})
\textit{Let $G$ be a $p$-group of order $p^n$. Then $s_k(G)\equiv 1 \ (mod \ p),$  $\forall \ k\in \lbrace 1,2,\ldots, n-1\rbrace.$}\\

\textbf{Theorem 2.1.2.} (see Theorem 1 of \cite{15})
\textit{Let $G$ be a $p$-group of order $p^n$ such that $G\not\cong\mathbb{Z}_{p^n}$ and $p$ is an odd prime number. Then $s_k(G)\equiv 1+p \ (mod \ p^2)$, $\forall \ k\in \lbrace 1,2,\ldots,n-1\rbrace$.}\\

Since $|L(G)|=2+\sum\limits_{k=1}^{n-1}s_k(G)$, by the above results it follows that:
\begin{itemize}
\item[\textit{i)}] $|L(G)|\equiv n+1 \ (mod \ p)$; 
\item[\textit{ii)}] $|L(G)|\equiv (p+1)(n+1)+2 \ (mod \ p^2);$ note that this congruence holds only for finite $p$-groups $G$ satisfying the hypotheses of Theorem 2.1.2.
\end{itemize}

There is more to say if we work with finite abelian $p$-groups. First of all, it is well known that an abelian $p$-group $G$ of order $p^n$ is isomorphic to a direct product of cyclic $p$-groups, i.e.
$$G\cong \mathbb{Z}_{p^{d_1}}\times \mathbb{Z}_{p^{d_2}}\times\ldots\times\mathbb{Z}_{p^{d_k}},$$
where $d_1,d_2,\ldots,d_k$ are positive integers such that $1\leq d_1\leq d_2\leq\ldots\leq d_k$ and $d_1+d_2+\ldots+d_k=n$. In other words $d=(d_k,d_{k-1},\ldots,d_1,0,\ldots)$ is a partition of $n$ with $d_i\geq 1, \ \forall \ i\in\lbrace 1,2,\ldots,k\rbrace$, and $G$ is a finite abelian $p$-group of type $d$. An interesting and difficult problem is to count the subgroups of $G$ and an answer is given by several papers like \cite{4,5,8,11,29}. Here, we will only recall Lemma 1.4.1 of \cite{8}. Let $l=(l_k,l_{k-1},\ldots, l_1,0,\ldots)$ be a partition of a non-negative integer $m\leq n$, such that $0\leq l_i\leq d_i, \ \forall \ i\in\lbrace 1,2,\ldots, k\rbrace$. Then $H_l\cong \mathbb{Z}_{p^{l_1}}\times \mathbb{Z}_{p^{l_2}}\times\ldots\times\mathbb{Z}_{p^{l_k}}$ is a subgroup of $G$ of type $l$. The number of subgroups of $G$ that are isomorphic to $H_l$ is given by
\begin{align}\label{r1}
\prod\limits_{i=1}^{d_k}p^{l'_{k-i}(d'_{k-i+1}-l'_{k-i+1})}\bigg[\substack{d'_{k-i+1}-l'_{k-i} \\ l'_{k-i+1}-l'_{k-i}}\bigg]_p,
\end{align}
where $(l'_k, l'_{k-1},\ldots, l'_1,0,\ldots)$ and $(d'_k, d'_{k-1},\ldots, d'_1, 0 \ldots)$ are the conjugate partitions (with respect to Ferrers diagrams) of $l$ and $d$, respectively,  while  $\big[\substack{n \\ k}\big]_p$ is the Gaussian binomial coefficient given by
$$\big[\substack{n \\ k}\big]_p=\frac{\prod\limits_{i=1}^n(p^i-1)}{\prod\limits_{i=1}^k(p^i-1)\prod\limits_{i=1}^{n-k}(p^i-1)}.$$
It is clear that it is difficult to work with formula \eqref{r1}, especially if one is interested in finding the total number of subgroups of $G$. Still, using different counting arguments, some explicit formulas that allow us to obtain the quantity $|L(G)|$ were given for abelian $p$-groups of rank 2 and 3. For more details, the reader may consult \cite{9,14,18,24,26}. Here, we only recall that the number of subgroups of $G\cong\mathbb{Z}_{p^{d_1}}\times\mathbb{Z}_{p^{d_2}}$, where $1\leq d_1\leq d_2$, may be computed using the explicit formula
\begin{align}\label{r6}
|L(G)|=\frac{1}{(p-1)^2}[(d_2-d_1+1)p^{d_1+2}-(d_2-d_1-1)p^{d_1+1}-(d_1+d_2+3)p+(d_1+d_2+1)],
\end{align}
while the number of subgroups of $G\cong\mathbb{Z}_{p^{d_1}}\times\mathbb{Z}_{p^{d_2}}\times
\mathbb{Z}_{p^{d_3}}$, where $1\leq d_1\leq d_2\leq d_3$, is given by 
\begin{align}\label{r7}
L(G)=\frac{A}{(p^2-1)^2(p-1)},
\end{align}
where
\begin{align*}
A=&(d_1+1)(d_3-d_2+1)p^{d_1+d_2+5}+2(d_1+1)p^{d_1+d_2+4}-2(d_1+1)(d_3-d_2)p^{d_1+d_2+3}-\\ &
2(d_1+1)p^{d_1+d_2+2}+(d_1+1)(d_3-d_2-1)p^{d_1+d_2+1}-(d_3+d_2-d_1+3)p^{2d_1+4}-2p^{2d_1+3}+ \\ &
(d_3+d_2-d_1-1)p^{2d_1+2}+(d_1+d_2+d_3+5)p^2+2p-(d_1+d_2+d_3+1).
\end{align*}

A relevant consequence of (\ref{r1}) is that for all $k\in\lbrace 1,2,\ldots,n-1\rbrace$, $s_k(G)$ may be viewed as a polynomial in $p$ with non-negative integer coefficients. Moreover, in his unpublished work, P. Hall proved that the number of subgroups of $G$ that are isomorphic to $H_l$ is equal to the number of subgroups of $G$ that are isomorphic to $\frac{G}{H_l}$. This result is also presented on p. 188 of \cite{16} and, as a consequence, for abelian $p$-groups $G$ or order $p^n$ we have
\begin{align}\label{r2}
s_k(G)=s_{n-k}(G), \ \forall \ k\in\bigg\lbrace 1,2,\ldots,\bigg[\frac{n}{2}\bigg]\bigg\rbrace.
\end{align}    

Finally, we recall the main theorem of \cite{7} which is also related to the quantities $s_k(G)$ of an  abelian $p$-group $G$ of order $p^n$.\\

\textbf{Theorem 2.1.3.} \textit{Let $G$ be an abelian $p$-group of order $p^n$. Then the sequence $(s_k(G))_{k=\overline{0,[\frac{n}{2}]}}$ is unimodal.}\\

In other words, the last result states that the finite sequence formed of the polynomials $s_k(G)$, where $k\in\lbrace 0,1,2,\ldots, \big[\frac{n}{2}\big]\rbrace$, has the following property:
$$\forall \ k\in \bigg\lbrace 1,2,\ldots,\bigg[\frac{n}{2}\bigg]\bigg\rbrace, s_k(G)-s_{k-1}(G) \ \text{is a polynomial in $p$ with non-negative integer coefficients.}$$
   
For a finite abelian $p$-group $G$ or order $p^n$ such that $G\not\cong \mathbb{Z}_{p^n}$, there are two ways to express the quantity $|L(G)|$ based on the parity of $n$. More exactly, if $p$ is any prime number, according to Theorem 2.1.1, for each $k\in\lbrace 1,2,\ldots,n-1\rbrace$, there is a non-negative integer $M_k$ such that $s_k(G)=M_kp+1.$ Moreover, since $G$ is a finite non-cyclic abelian $p$-group, we can take $M_k\geq 1,$ $\ \forall \ k\in\lbrace 1,2,\ldots,n-1\rbrace$, as a consequence of Proposition 1.3 of \cite{3}. Then, using (\ref{r2}), we have $M_1=M_{n-1}, M_2=M_{n-2},\ldots, M_{[\frac{n}{2}]}=M_{n-[\frac{n}{2}]}$. Therefore, the number of subgroups of $G$ is given by 
\begin{align}\label{r3}
|L(G)|=\begin{cases} 2(M_1+M_2+\ldots+M_{\frac{n-1}{2}})p+(n+1) &\mbox{, if } n\equiv 1 \ (mod \ 2) \\  [2(M_1+M_2+\ldots+M_{\frac{n}{2}-1})+M_{\frac{n}{2}}]p+(n+1) &\mbox{, if } n\equiv 0 \ (mod \ 2) \end{cases}.
\end{align}
Similarly, if $G$ is a $p$-group of order $p^n$ such that $G\not\cong\mathbb{Z}_{p^n}$ and $p$ is odd, by Theorem 2.1.2, it follows that for each $k\in\lbrace 1,2,\ldots,n-1\rbrace$, there is a non-negative integer $N_k$ such that $s_k(G)=N_kp^2+p+1$. Following the same reasoning as the one that was done to obtain (\ref{r3}) , we have
\begin{align}\label{r4}
|L(G)|=\begin{cases} 2(N_1+N_2+\ldots+N_{\frac{n-1}{2}})p^2+(n-1)(p+1)+2 &\mbox{, if } n\equiv 1 \ (mod \ 2) \\  [2(N_1+N_2+\ldots+N_{\frac{n}{2}-1})+N_{\frac{n}{2}}]p^2+(n-1)(p+1)+2 &\mbox{, if } n\equiv 0 \ (mod \ 2) \end{cases}.
\end{align}
In the end, we remark that $1\leq M_1\leq M_2\leq\ldots\leq M_{[\frac{n}{2}]}$ and $0\leq N_1\leq N_2\leq\ldots\leq N_{[\frac{n}{2}]}$, as Theorem 2.1.3 indicates.   

\subsection{Basic properties of $\beta$}

Any finite group $G$ has at least 2 subgroups, and, by Corollary 1.6 of \cite{6}, we have 
$$|L(G)|\leq |G|^{(\frac{1}{4}+o(1))log_2|G|},$$  
this upper bound for the number of subgroups of $G$ being the best possible one, in the sense that it is ``close" to $|L(\mathbb{Z}_2^n)|$. It follows that 
$$\frac{2}{|G|}\leq \beta(G)\leq \frac{|G|^{(\frac{1}{4}+o(1))log_2|G|}}{|G|}.$$
Once that $|G|$ increases, the lower bound goes to 0, while the upper one approaches infinity. Hence, it is clear that 
$$\beta(G)\in (0,\infty).$$

If two finite groups $G_1$ and $G_2$ are isomorphic, then $\beta(G_1)=\beta(G_2)$. The converse is false since 
$$\beta(\mathbb{Z}_2\times\mathbb{Z}_8)=\beta(M_{16})=\beta(Q_{16})=\frac{11}{16},$$
and any two of the above three groups are not isomorphic.

An important property for our study is the multiplicativity of $\beta$. This property plays a significant role when it comes to prove some density results related to $\beta$ in Section 4.\\

\textbf{Proposition 2.2.1.} 
\textit{Let $(G_i)_{i=\overline{1,k}}$ be a family of finite groups having coprime orders. Then
$$\beta(\xmare{i=1}k G_i)=\prod\limits_{i=1}^k\beta(G_i).$$}

A direct consequence of this multiplicativity is expressing the quantity $\beta$ associated to any finite nilpotent group $G$ since $G$ can be written as the direct product of its Sylow subgroups. So, if $(P_i)_{i=\overline{1,k}}$ are the Sylow subgroups of $G$, we have
$$\beta(G)=\prod\limits_{i=1}^k\beta(P_i).$$
  
One may prove that any finite group $G$ having at most 5 subgroups is abelian. There are different ways to check this, but we give a proof which is connected with some of the results that were recalled in the previous Subsection. Also, we write the above property in terms of $\beta$.\\

\textbf{Proposition 2.2.2.} 
\textit{Let $G$ be a finite group. If $\beta(G)\leq \frac{5}{|G|}$, then $G$ is abelian.}

\textbf{Proof.}
Suppose that $G$ is non-abelian. If $G$ is a $p$-group of order $p^n$, since $G$ would have at least one $p$-subgroup of order $p^k$ for each $k\leq n$ and $G$ is non-abelian, it follows that $n\in\lbrace
3,4\rbrace$. Then $|L(G)|\geq 3+s_1(G)$. By Theorem 2.1.1, we have $s_1(G)\geq p+1$ as a consequence of the fact that $G$ is non-cyclic. Then $\beta(G)\geq \frac{p+4}{|G|}\geq \frac{6}{|G|},$ a contradiction. 
  
If $G$ is not a $p$-group, there are exactly 2 distinct prime numbers, say $p$ and $q$, that are divisors of $|G|$. Indeed, if at least 4 such divisors exist, then $\beta(G)>\frac{5}{|G|}$ as a consequence of Cauchy's theorem, a contradiction. If there are 3 distinct prime divisors of $|G|$, then $G$ has 3 Sylow subgroups of different orders and all of them must be normal. Otherwise Sylow's 3rd theorem would imply that $\beta(G)>\frac{5}{|G|}$, a contradiction. But, since all normal subgroups are permutable, again we would arrive at the same contradiction. Therefore, $G$ has 2 Sylow subgroups $P$ and $Q$ of orders $p^x$ and $q^y$, respectively, where $x$ and $y$ are positive integers. Then 
$\beta(G)\geq \frac{x+y+2}{|G|},$ since $P$ and $Q$ have at least $x$ and $y$ non-trivial subgroups, respectively. If $x\geq 3$ or $y\geq 3$, we contradict our hypothesis. The cases $(x,y)\in \lbrace (1,2), (2,1)\rbrace$ are excluded using the same reasoning involving Sylow's 3rd theorem and the fact that normality implies permutability. Then $x=y=1$, so $G$ is a group of order $pq$. Since $G$ is non-abelian, it follows that $P$ or $Q$ is not a normal subgroup of $G$. Then, once again $\beta(G)>\frac{5}{|G|}$ as a consequence of Sylow's 3rd theorem, a contradiction. Hence, our assumption that $G$ is non-abelian is false and the proof is finished. 
\hfill\rule{1,5mm}{1,5mm}
 
\section{Bounds and minima problems related to $p$-groups}

The aim of this Section is to add some properties of $\beta$ by restricting our study to finite $p$-groups. We start by indicating some bounds for this quantity and we continue by solving some minima problems, the main result pointing out the finite $p$-groups for which $\beta$ attains its second minimum value.\\   

We denote by $\mathcal{P}$ the class of $p$-groups of order $p^n$, where $n\geq 3$. Let $G\in\mathcal{P}$. Then
$$\beta(G)=\frac{\sum\limits_{k=0}^n s_k(G)}{p^n}.$$
Since $G$ has at least one subgroup of order $p^k$, $k\in\lbrace 1,2,\ldots,n-1\rbrace$, the minimum value of $\beta$, on $\mathcal{P}$, is attained if and only if $s_k(G)=1$, $ \ \forall \ k\in\lbrace 1,2,\ldots, n-1\rbrace$, i.e., if and only if $G\cong \mathbb{Z}_{p^n}$. Using Theorem 5.17 of \cite{3}, we deduce that $\beta$ attains its maximum value on $\mathcal{P}$ if and only if $G\cong\mathbb{Z}_p^n.$ Hence, for all $G\in\mathcal{P}$, we have
$$\beta(\mathbb{Z}_{p^n})\leq \beta(G)\leq \beta(\mathbb{Z}_p^n).$$

A connection between $G\in\mathcal{P}$ and some of its quotients may be written as a consequence of Theorem 1.3 of \cite{20}. More exactly, let $G\in\mathcal{P}$ and $H$ be a normal subgroup of $G$ such that $|H|=p$. Then,
$$\beta(G)\leq \beta\bigg(\frac{G}{H}\times\mathbb{Z}_p\bigg).$$
In the same paper (see Theorem 1.4), the author proves that if $G$ is a non-elementary abelian $p$-group of order $p^n$, where $p$ is odd and $n\geq 3$ is a positive integer, then 
$$s_k(G)\leq s_k(M_p(1,1,1)\times \mathbb{Z}_p^{n-3}), \ \forall \ k\in\lbrace 1,2,\ldots, n-1\rbrace,$$
where
$$M_p(1,1,1)=\langle a,b,c \ | \ a^p=b^p=c^p=1, [a,b]=c, [c,a]=[c,b]=1\rangle.$$
He conjectures that this result also holds for $p=2$ and a proof of this fact is given in \cite{27}. This means that the second maximum value of $\beta$, on $\mathcal{P}$, is attained when one works with the group $M_p(1,1,1)\times \mathbb{Z}_p^{n-3}.$\\

\textit{What about the second minimum value?} In what follows, we provide an answer to this question. In this regard, it is worth to recall Theorem 2.2 of \cite{19}. This result states that for a finite $p$-group $G$ of order $p^n$, we have $s_k(G)=p+1, \ \forall \ k\in\lbrace 1,2,\ldots,n-1\rbrace$ if and only if $G\cong\mathbb{Z}_p\times\mathbb{Z}_{p^{n-1}}$ or $G\cong M(p^n)$.\\ 

Firstly, we indicate an answer to the above question if we work only with abelian groups contained in $\mathcal{P}$.\\

\textbf{Proposition 3.1.} \textit{Let $G\in\mathcal{P}$ such that $G$ is abelian and $G\not\cong\mathbb{Z}_{p^n}$. Then $$\beta(G)\geq \beta(\mathbb{Z}_p\times\mathbb{Z}_{p^{n-1}}).$$
The equality holds if and only if $G\cong\mathbb{Z}_p\times\mathbb{Z}_{p^{n-1}}.$}   

\textbf{Proof.} Let $G$ be a group as indicated by our hypothesis and suppose that $\beta(G)< \beta(\mathbb{Z}_p\times\mathbb{Z}_{p^{n-1}})$. By (\ref{r6}), we have
$$\beta(\mathbb{Z}_p\times\mathbb{Z}_{p^{n-1}})=\frac{(n-1)(p+1)+2}{p^n}.$$ 
If $n$ is odd, then according to (\ref{r3}), we deduce that
\begin{align*}
\beta(G)<\beta(\mathbb{Z}_p\times\mathbb{Z}_{p^{n-1}}) \Longrightarrow 2(M_1+M_2+\ldots+M_{\frac{n-1}{2}})< n-1,
\end{align*}
where $1\leq M_1\leq M_2\leq\ldots\leq M_{\frac{n-1}{2}}$. But, 
$$2(M_1+M_2+\ldots+M_{\frac{n-1}{2}})\geq 2\cdot\frac{n-1}{2}=n-1,$$ 
and this leads to a contradiction. Consequently, the inequality $\beta(G)\geq \beta(\mathbb{Z}_p\times\mathbb{Z}_{p^{n-1}})$ holds. Following the above reasoning, one can analyse the case of even integers $n\geq 3$ and arrive at the same conclusion.

Concerning the situation where equality holds, we have 
$$\beta(G)=\beta(\mathbb{Z}_p\times\mathbb{Z}_{p^{n-1}})\Longleftrightarrow M_1+M_2+\ldots+M_{n-1}=n-1,$$
where $M_k\geq 1, \ \forall \ k\in\lbrace 1,2,\ldots,n-1\rbrace$. Then, 
\begin{align*}
\beta(G)=\beta(\mathbb{Z}_p\times\mathbb{Z}_{p^{n-1}})&\Longleftrightarrow M_k=1, \ \forall \ k\in\lbrace 1,2,\ldots,n-1\rbrace\\&\Longleftrightarrow s_k(G)=p+1, \ \forall \ k\in\lbrace 1,2,\ldots,n-1\rbrace\\&\Longleftrightarrow G\cong\mathbb{Z}_p\times\mathbb{Z}_{p^{n-1}}
\end{align*}      
We note that the last equivalence is a consequence of the classification that was recalled above.     
\hfill\rule{1,5mm}{1,5mm}\\

In what follows, we find the second minimum of $\beta$ on the entire $\mathcal{P}$. It is quite interesting that, in some cases, there are at least 2 minimum points in $\mathcal{P}$ associated to this minima problem.\\

\textbf{Theorem 3.2.} \textit{Let $G\in\mathcal{P}$ such that $G\not\cong\mathbb{Z}_{p^n}$.
\begin{itemize}
\item[i)] If $p$ is odd, then $\beta(G)\geq\beta(\mathbb{Z}_p\times\mathbb{Z}_{p^{n-1}})=\beta(M(p^n)).$
\item[ii)] If $p=2$, then 
$\begin{cases} \beta(G)\geq \beta(Q_8) &\mbox{, if } n=3 \\ \beta(G)\geq \beta(\mathbb{Z}_2\times\mathbb{Z}_8)=\beta(Q_{16})=\beta(M(16)) &\mbox{, if } n=4 \\ \beta(G)\geq \beta(\mathbb{Z}_2\times\mathbb{Z}_{2^{n-1}})=\beta(M(2^n)) &\mbox{, if } n>4\end{cases}.$
\end{itemize}
The equality holds if and only if $G$ is isomorphic to one of the indicated minimum points corresponding to each case.}

\textbf{Proof.}
Let $G\in\mathcal{P}$ with $G\not\cong\mathbb{Z}_{p^n}$. We must find another group $G_1$ having the same properties such that $\beta(G)\geq\beta(G_1)$. To obtain the smallest value of $\beta(G_1)$, we must lower the quantities $s_k(G_1), \ \forall \ k\in\lbrace 1,2,\ldots, n-1\rbrace.$ If there is an integer $k$ such that $s_k(G_1)=1$, since $G_1$ is non-cyclic, by Proposition 1.3 of \cite{3}, the previous equality holds if and only if $G_1\cong Q_{2^n}$. Otherwise, if $s_k(G_1)>1, \ \forall \ k\in\lbrace 1,2,\ldots,n-1\rbrace$, then Theorem 2.1.1 indicates that the lowest possible value of $s_k(G_1)$ is $p+1$ for all $k\in\lbrace 1,2,\ldots,n-1\rbrace$. As we previously remarked, this happens if and only if $G_1\cong\mathbb{Z}_p\times\mathbb{Z}_{p^{n-1}}$ or $G_1\cong M(p^n)$.

Hence, if $p$ is odd, then the minimum point $G_1$ is isomorphic to $\mathbb{Z}_p\times\mathbb{Z}_{p^{n-1}}$ or $M(p^n)$ and the second minimum value of $\beta$, on $\mathcal{P}$, is
$$\beta(\mathbb{Z}_p\times\mathbb{Z}_{p^{n-1}})=\beta(M(p^n))=\frac{(n-1)(p+1)+2}{p^n}.$$

If $p=2$, the determination of the minimum point $G_1$ is related to the value of $n$. We recall that 
$$|L(Q_{2^n})|=2^{n-1}+n-1,$$ this result being indicated in \cite{23}. We have 
$$\beta(\mathbb{Z}_2\times\mathbb{Z}_{2^{n-1}})=\beta(M(2^n))=\frac{3(n-1)+2}{2^n} \text{ \ and \ } \beta(Q_{2^n})=\frac{2^{n-1}+n-1}{2^n}.$$
It is easy to check that $\beta(Q_{2^n})\geq\beta(\mathbb{Z}_2\times\mathbb{Z}_{2^{n-1}})=\beta(M(2^n))$ for any integer $n\geq 4$ and that the equality holds for $n=4$. Then, in order to finish the proof, we distinguish the following 3 cases:
\begin{itemize}
\item[--] If $n=3$, then the minimum point is $G_1\cong Q_8$ and  
$$\beta(G)\geq\beta(Q_8)=\frac{3}{4}.$$
\item[--] If $n=4$, then the minimum point $G_1$ is isomorphic to $\mathbb{Z}_2\times\mathbb{Z}_8, Q_{16}$ or $M(16)$ and  
$$\beta(G)\geq\beta(\mathbb{Z}_2\times\mathbb{Z}_8)=\beta(Q_{16})=\beta(M(16))=\frac{11}{16}.$$
\item[--] If $n>4$, then the minimum point $G_1$ is isomorphic to $\mathbb{Z}_2\times\mathbb{Z}_{2^{n-1}}$ or $M(2^n)$ and 
$$\beta(G)\geq\beta(\mathbb{Z}_2\times\mathbb{Z}_{2^{n-1}})=\beta(M(2^n))=\frac{3(n-1)+2}{2^n}.$$
\end{itemize}
Since we determined all possible minimum points $G_1$ corresponding to each case, the equality $\beta(G)=\beta(G_1)$ holds if and only if $G\cong G_1$.
\hfill\rule{1,5mm}{1,5mm}\\

The results that were proved in this Section may be also interpreted as follows.\\

\textbf{Corrolary 3.3.} \textit{Let $G\in\mathcal{P}$ such that $G$ is abelian. Then
$$\beta(G)\leq\frac{(n-1)(p+1)+2}{p^n}\Longleftrightarrow G\cong \mathbb{Z}_{p^n} \text{ or } G\cong\mathbb{Z}_p\times\mathbb{Z}_{p^{n-1}}.$$} 

\textbf{Corrolary 3.4.} \textit{Let $G\in\mathcal{P}$.
\begin{itemize}
\item[i)] If $p$ is odd, then 
$$\beta(G)\leq\frac{(n-1)(p+1)+2}{p^n}\Longleftrightarrow G\cong\mathbb{Z}_{p^n}, G\cong\mathbb{Z}_p\times\mathbb{Z}_{p^{n-1}} \text{ or } G\cong M(p^n).$$
\item[ii)] If $p=2$, then 
$\begin{cases} \beta(G)\leq\frac{3}{4}\Longleftrightarrow G\cong\mathbb{Z}_8 \text{ or } G\cong Q_8 &\mbox{, for } n=3 \\ \beta(G)\leq\frac{11}{16} \Longleftrightarrow G\cong\mathbb{Z}_{16}, G\cong\mathbb{Z}_2\times\mathbb{Z}_8, G\cong Q_{16} \text{ or } G\cong M(16) &\mbox{, for } n=4 \\ \beta(G)\leq \frac{3(n-1)+2}{2^n}\Longleftrightarrow G\cong\mathbb{Z}_{2^n}, G\cong\mathbb{Z}_2\times\mathbb{Z}_{2^{n-1}} \text{ or } G\cong M(2^n) &\mbox{, for } n>4\end{cases}.$
\end{itemize}}

One may go further and try to find the third minimum (maximum) value of $\beta$, on the class of $p$-groups of order $p^n$, where $n\geq 4$. In this regard, we prove a result which may be considered as a starting point for such a study. More exactly, we show that the third minimum value of $\beta$, on the class of abelian $p$-groups of order $p^n$ such that $n\geq 4$ and $p$ is odd, is attained at the ``point" $\mathbb{Z}_{p^2}\times\mathbb{Z}_{p^{n-2}}$. \\

We recall that for a positive integer $n$, there is a bijection between the set of partitions of $n$ and the set of types of abelian $p$-groups of order $p^n$. More exactly, for a partition $d=(d_k,d_{k-1},\ldots,d_1,0,\ldots)$ of $n$, there is an unique abelian $p$-group $G\cong\mathbb{Z}_{p^{d_1}}\times\mathbb{Z}_{p^{d_2}}\times\ldots\times\mathbb{Z}_{p^{d_k}}$ of type $d$ and order $p^n$. Moreover, the relation ``$\preceq$" defined by
\begin{align*}
(d_k,d_{k-1},\ldots,d_1,0,\ldots)\preceq (e_k,e_{k-1},\ldots,e_1,0,\ldots)\Longleftrightarrow & d_i=e_i, \ \forall \ i\in\lbrace 1,2,\ldots,k\rbrace\\ & \text{ or } \\ & \exists \ i\in\lbrace 1,2,\ldots,k-1\rbrace \text{ such that } d_i<e_i \\ & \text{ and } d_{i+1}=e_{i+1},\ldots,d_k=e_k,
\end{align*}
is a total order on the set containing all partitions of $n$.\\   

We are ready to prove a preliminary result that will also be helpful in Section 5.\\

\textbf{Lemma 3.5.} \textit{$\beta$ is strictly decreasing on the class of abelian $p$-groups of rank 2 and order $p^n$, where $n\geq 4$.}

\textbf{Proof.} Let $n\geq 4$ be a positive integer. Without loss of generality, we choose the abelian $p$-group $G_1\cong\mathbb{Z}_{p^{d_1}}\times\mathbb{Z}_{p^{d_2}}$ of type $d=(d_2,d_1,0,\ldots)$ and order $p^n$, where $2\leq d_1\leq d_2$. Since the set of partitions of $n$ is totally ordered, it is sufficient to take the consecutive partition of $d$ with respect to ``$\preceq$", i.e. $(d_2+1,d_1-1,0,\ldots)$, and the corresponding group $G_2\cong\mathbb{Z}_{p^{d_1-1}}\times\mathbb{Z}_{p^{d_2+1}}$, and prove that $\beta(G_2)<\beta(G_1)$. Using (\ref{r6}), we have
\begin{align*}
\beta(G_2)<\beta(G_1)&\Longleftrightarrow (d_2-d_1+3)p-(d_2-d_1+1)<(d_2-d_1+1)p^2-(d_2-d_1-1)p\\
& \Longleftrightarrow 2p<p^2+1. 
\end{align*}
Since the last inequality holds for any prime $p$, the proof is complete. 
\hfill\rule{1,5mm}{1,5mm}\\
 
We remark that, in general, $\beta$ is not monotonic on the class of finite abelian $p$-groups of a given order. For instance, if $p=2$ and $n=9$, then we have
$(4,4,1,0,\ldots)\preceq (5,2,2,0,\ldots)\preceq (5,3,1,0,\ldots),$
but $\beta(\mathbb{Z}_2\times\mathbb{Z}_{16}\times\mathbb{Z}_{16})=322, \beta(\mathbb{Z}_4\times\mathbb{Z}_4\times\mathbb{Z}_{32})=354$ and $\beta(\mathbb{Z}_2\times\mathbb{Z}_8\times\mathbb{Z}_{32})=258.$ We mention that these numbers may be obtained using (\ref{r7}) or GAP \cite{28}.\\
   
Finally, before proving the last result of this Section, we mention that for any finite $p$-group $G$ of order $p^n$, there is a bijection between the set containing all maximal subgroups of $G$ and the set formed of the maximal subgroups of $\frac{G}{\Phi(G)}\cong\mathbb{Z}_p^{d(G)}$. Here, we denoted the Frattini subgroup and the minimal number of generators of $G$ by $\Phi(G)$ and $d(G)$, respectively. Hence, 
$$s_{n-1}(G)=s_{d(G)-1}(\mathbb{Z}_p^{d(G)})=\frac{p^{d(G)}-1}{p-1}.$$     

\textbf{Proposition 3.5.} \textit{Let $G$ be an abelian $p$-group of order $p^n$ such that $n\geq 4$ and $p$ is odd. Suppose that $G\not\cong\mathbb{Z}_{p^n}$ and $G\not\cong\mathbb{Z}_p\times\mathbb{Z}_{p^{n-1}}.$ Then
$$\beta(G)\geq\beta(\mathbb{Z}_{p^2}\times\mathbb{Z}_{p^{n-2}}).$$}

\textbf{Proof.} Let $G$ be a group as indicated by our hypothesis. Suppose that $\beta(G)<\beta(\mathbb{Z}_{p^2}\times\mathbb{Z}_{p^{n-2}}).$ We choose $n\geq 4$ to be an even integer and we mention that one may follow the same reasoning to complete the proof for odd integers. By (\ref{r6}), we have 
$$\beta(\mathbb{Z}_{p^2}\times\mathbb{Z}_{p^{n-2}})=\frac{(n-3)p^2+(n-1)(p+1)+2}{p^n}.$$
Therefore, since $G$ is non-cyclic and $p$ is odd, we can use (\ref{r4}) to deduce that
$$\beta(G)<\beta(\mathbb{Z}_{p^2}\times\mathbb{Z}_{p^{n-2}})\Longleftrightarrow 2(N_1+N_2+\ldots+N_{\frac{n}{2}-1})+N_{\frac{n}{2}}<n-3,$$
where $0\leq N_1\leq N_2\leq\ldots\leq N_{\frac{n}{2}}.$ If $N_1\geq 1$, then
$$2(N_1+N_2+\ldots+N_{\frac{n}{2}-1})+N_{\frac{n}{2}}\geq n-1,$$
a contradiction. Then $N_1=0$, so, by (\ref{r2}), we have $s_1(G)=s_{n-1}(G)=p+1$. But, as we mentioned above, we also have $s_{n-1}(G)=\frac{p^{d(G)}-1}{p-1}$. It follows that $d(G)=2$. Consequently, $G\cong\mathbb{Z}_{p^{d_1}}\times\mathbb{Z}_{p^{d_2}}$, where $1\leq d_1\leq d_2.$ According to Lemma 3.4, since $\beta(G)<\beta(\mathbb{Z}_{p^2}\times\mathbb{Z}_{p^{n-2}})$, it follows that $d_1=1$ and $d_2=n-1$. Then $G\cong\mathbb{Z}_p\times\mathbb{Z}_{p^{n-1}}$ and we contradict our hypothesis. Therefore $\beta(G)\geq\beta(\mathbb{Z}_{p^2}\times\mathbb{Z}_{p^{n-2}}),$ as desired.   
\hfill\rule{1,5mm}{1,5mm}\\

\section{Some density results associated to $\beta$}

We denote by $\mathcal{F}$ the class of all finite groups. Let $\mathcal{A}$ be the subclass of $\mathcal{F}$ containing all finite abelian groups. In this Section, our main aim is to prove that the sets
$$\lbrace \beta(G) \ | \ G\in\mathcal{A}\rbrace \text{ and } \lbrace \beta(G) \ | \ G\in\mathcal{F}\rbrace$$
are dense in $[0,\infty).$\\

A first step to reach our purpose is to recall the Proposition marked on p. 863 of \cite{17}. This result mainly states that given a sequence $(x_n)_{n\geq 1}$ of positive real numbers such that $\displaystyle \lim_{n \to\infty}x_n=0$ and $\sum\limits_{n=1}^{\infty}x_n$ is divergent, then $\sum(\lbrace x_i\rbrace_{i=1}^\infty)=[0,\infty)$. As the author indicates, by $\sum(\lbrace x_i\rbrace_{i=1}^\infty)$ we denote the set containing the sums of all finite and summable infinite subsequences of $x_n$, as well as the sum of the empty subsequence of $x_n$, which is 0.\\

Taking into consideration the above statements, we are ready to prove the following preliminary result.\\

\textbf{Lemma 4.1.} \textit{Let $(x_n)_{n\geq 1}$ be a sequence of positive real numbers such that $\displaystyle \lim_{n \to\infty}x_n=0$ and $\sum\limits_{n=1}^{\infty}x_n$ is divergent. Then the set containing the sums of all finite subsequences of $(x_n)$ is dense in $[0,\infty)$.}

\textbf{Proof.} Let $x\in[0,\infty)$ and let $(x_n)_{n\geq 1}$ be a sequence that satisfies our hypothesis. According to the result that we marked before we started this proof, we know that there is a subsequence $(x_{n_k})_{k\geq 1}$ of $x_n$ such that the sum of its elements is equal to $x$. If $x_{n_k}$ is finite, then the conclusion follows since we may take the constant sequence formed of the sums $\sum x_{n_k}$ which is convergent to $x$. If $x_{n_k}$ is summable infinite, then we know that $\sum\limits_{k=1}^{\infty}x_{n_k}=x$. This may be written as $\displaystyle \lim_{N \to\infty}\sum\limits_{k=1}^N x_{n_k}=x$. Then we may choose the sequence, indexed by $N$, formed by the sums $\sum\limits_{k=1}^N x_{n_k}$, where $N$ is sufficiently large, which again is convergent to $x$. Since the subsequences $(x_{n_k})_{k=\overline{1,N}}$ are finite, the proof is complete.
\hfill\rule{1,5mm}{1,5mm}\\

Before proving another density result, we recall that the series $\sum\limits_{p=prime}\frac{1}{p}$ is divergent. Also, for a continuous function $f:\mathbb{R}\longrightarrow\mathbb{R}$ and two subsets $A$ and $B$ of $\mathbb{R}$ such that $\overline{A}=\overline{B}$, we have $\overline{f(A)}=\overline{f(B)}.$ \\

\textbf{Proposition 4.2.} \textit{The set $$\lbrace \beta(\xmare{i\in I}{ }\mathbb{Z}_{p_i}^4) \ | \ I\subset \mathbb{N}\setminus\lbrace 0\rbrace, |I|<\infty \text{ and } p_i \text{ is the ith prime number}, \ \forall \ i\in I\rbrace$$ is dense in $[1,\infty)$.}

\textbf{Proof.} Let $p_n$ be the $n$th prime number, where $n\geq 1$ is a positive integer. We have
$$\beta(\mathbb{Z}_{p_n}^4)=\sum\limits_{k=0}^4\big[\substack{4 \\ k}\big]_{p_n}=\frac{p_n^4+3p_n^3+4p_n^2+3p_n+5}{p_n^4}.$$
We check that the sequence $(x_n)_{n\geq 1}$, where $x_n=\ln\beta(\mathbb{Z}_{p_n}^4), \ \forall \ n\geq 1$, satisfies the hypotheses of Lemma 4.1. It is clear that $x_n>0, \ \forall \ n\geq 1.$ Also, we have
$\displaystyle \lim_{n \to\infty}x_n=\displaystyle \lim_{n \to\infty}\ln\beta(\mathbb{Z}_{p_n}^4)=\ln 1=0.$
Finally, the function $f:[2,\infty)\longrightarrow\mathbb{R}$ given by $f(x)=\ln\bigg(\frac{x^4+3x^3+4x^2+3x+5}{x^4}\bigg)-\frac{1}{x}$ takes positive values since 
$f'(x)=-\frac{2x^4+5x^3+5x^2+17x-5}{x^2(x^4+3x^3+4x^2+3x+5)}<0, \ \forall \ x\geq 2,$ and $\displaystyle \lim_{x \to\infty}f(x)=0$. Consequently, $f(p_n)>0, \ \forall \ n\geq 1$, so $x_n>\frac{1}{p_n}, \ \forall \ n\geq 1$. Therefore, the series $\sum\limits_{n=1}^\infty x_n$ is divergent.

By Lemma 4.1, we deduce that
$$\overline{\lbrace \ln\prod\limits_{i\in I}\beta(\mathbb{Z}_{p_i}^4) \ | \ I\subset \mathbb{N}\setminus\lbrace 0\rbrace, |I|<\infty \text{ and } p_i \text{ is the \textit{i}th prime number}, \ \forall \ i\in I\rbrace}=[0,\infty).$$  
Since the exponential function is continuous and the above relation expresses an equality between the closures of two sets of $\mathbb{R}$, we obtain
$$\overline{\lbrace \prod\limits_{i\in I}\beta(\mathbb{Z}_{p_i}^4) \ | \ I\subset \mathbb{N}\setminus\lbrace 0\rbrace, |I|<\infty \text{ and } p_i \text{ is the \textit{i}th prime number}, \ \forall \ i\in I\rbrace}=[1,\infty).$$
The conclusion follows by applying the multiplicativity of $\beta$ (Proposition 2.2.1) for the left-hand side of the above equality.   
\hfill\rule{1,5mm}{1,5mm}\\    

Let $k$ be a positive integer. Note that Lemma 4.1 also holds if we work with the sequence $(x_n)_{n\geq 1,n\ne k}$ because this sequence contains only positive real numbers, its limit is 0 and the series 
$\sum\limits_{n=1,n\ne k}^\infty x_n$ is divergent. This happens since the nature and the limit of $x_n$, as well as the nature of the series $\sum\limits_{n=1}^\infty x_n$, are not affected by eliminating a finite number of terms of $x_n$. Consequently, Proposition 4.2 also holds in this case and may be rewritten as 
$$\overline{\lbrace \beta(\xmare{i\in I}{ }\mathbb{Z}_{p_i}^4) \ | \ I\subset \mathbb{N}\setminus\lbrace 0,k\rbrace, |I|<\infty \text{ and } p_i \text{ is the ith prime number}, \ \forall \ i\in I\rbrace}=[1,\infty).$$
We have all necessary ingredients to prove the main result of this Section.\\

\textbf{Theorem 4.3.} \textit{\begin{itemize}
\item[i)] The set $\lbrace \beta(G) \ | \ G\in\mathcal{A}\rbrace$ is dense in $[0,\infty).$
\item[ii)] The set $\lbrace \beta(G) \ | \ G\in\mathcal{F}\rbrace$ is dense in $[0,\infty).$
\end{itemize}}

\textbf{Proof.} \textit{i)} Once again, denote by $p_n$ the $n$th prime number. According to Proposition 4.2, each $x\in [1,\infty)$ is an adherent point of the set $\lbrace \beta(G) \ | \ G\in\mathcal{A}\rbrace.$ Obviously, $0$ is also an adherent point of this set since $\displaystyle \lim_{n \to\infty}\beta(\mathbb{Z}_{p_n})=\lim_{n \to\infty}\frac{2}{p_n}=0.$

Let $x\in(0,1)$. We choose $p_k$ to be the first prime number for which the inequality $\frac{2}{p_k}\leq x$. Then $\frac{p_k}{2}\cdot x\in [1,\infty)$. So, according to the above rewriting of Proposition 4.2, there is a sequence of finite subsets $(I_n)_{n\geq 1}$ of $\mathbb{N}\setminus\lbrace 0,k\rbrace$, or equivalently, there is a sequence of finite abelian groups $(G_n)_{n\geq 1}$, where $G_n=\xmare{i\in I_n}{ }\mathbb{Z}_{p_i}^4$, such that $\displaystyle \lim_{n \to\infty}\beta(G_n)=\frac{p_k}{2}\cdot x.$ Finally, we consider the sequence $(\mathbb{Z}_{p_k}\times G_n)_{n\geq 1}$ and, by applying the multiplicativity of $\beta$, we obtain
$$\displaystyle \lim_{n \to\infty}\beta(\mathbb{Z}_{p_k}\times G_n)=\displaystyle \lim_{n \to\infty}\beta(\mathbb{Z}_{p_k})\beta(G_n)=\frac{2}{p_k}\cdot \frac{p_k}{2}\cdot x=x.$$
Hence $x$ is an adherent point of the set $\lbrace \beta(G) \ | \ G\in\mathcal{A}\rbrace.$ 

Therefore, $[0,\infty)\subseteq \overline{\lbrace \beta(G) \ | \ G\in\mathcal{A}\rbrace}$ and, since the converse inclusion is trivial, the conclusion follows.

\textit{ii)} We have $\lbrace \beta(G) \ | \ G\in\mathcal{A}\rbrace\subseteq \lbrace \beta(G) \ | \ G\in\mathcal{F}\rbrace\subseteq [0,\infty).$ Then, by taking the closures of these 3 sets, we deduce that the set $\lbrace \beta(G) \ | \ G\in\mathcal{F}\rbrace$ is dense in $[0,\infty).$ 
\hfill\rule{1,5mm}{1,5mm}
 
\section{$\beta$ viewed as a function on the subgroup lattice}

Until now, we studied some properties and problems related to $\beta$ viewed as a quantity associated to some remarkable classes of finite groups. Alternatively, given a finite group $G$, one may consider the function 
$$\beta:L(G)\longrightarrow (0,\infty) \text{ given by } H\mapsto \beta(H), \ \forall \ H\in L(G),$$
and study some properties connected with it.
 
For instance, in this Section, we work only with finite abelian $p$-groups and we provide a solution to the following problem:\\ 

\textit{Classify all finite abelian $p$-groups $G$ satisfying $\beta(H)\leq 1, \ \forall \ H\in L(G)$.}
\\
   
Concerning this question, one may ask why we compare the quantities $\beta(H)$, where $H$ is a subgroup of $G$, with 1. There are two reasons behind this choice. First of all, in the previous Section, we noted that $\beta(\mathbb{Z}_p^4)=\frac{p^4+3p^3+4p^2+3p+5}{p^4}$, and it is clear that this quantity is greater than 1 for any prime $p$. Then, any finite abelian $p$-group $G$ of rank $k\geq 4$ has a subgroup $H\cong\mathbb{Z}_p^4$ and we have $\beta(H)>1.$ Hence, to provide an anwer to our question, it is sufficient to study the finite abelian $p$-groups of rank $k\leq 3$. Secondly, by working at this classification, we also provide a partial solution to Problem 6.2 that is suggested in the following Section of the paper.\\

Once again we need to prove some preliminary results first. We note that minor computational details will be omitted.\\

\textbf{Lemma 5.1.} \textit{Let $G$ be a cyclic $p$-group or order $p^n$. Then $\beta(G)\leq 1$ and the equality holds if and only if $G\cong\mathbb{Z}_2$.}

\textbf{Proof.} Let $p$ be a prime number and let $G\cong \mathbb{Z}_{p^n}$, where $n$ is a positive integer. Then $\beta(G)=\frac{n+1}{p^n}$. The function $f:[1,\infty)\longrightarrow \mathbb{R}$ given by $f(x)=\frac{x+1}{p^x}$ is strictly decreasing since $f'(x)=\frac{1-(x+1)\ln p}{p^x}<0, \ \forall \ x\geq 1$. Then $f(x)\leq f(1), \ \forall \ x\in[1,\infty).$ Consequently, we remark that
$\beta(G)\leq\frac{2}{p}\leq 1$ and the equality $\beta(G)=1$ holds if and only if $n=1$ and $p=2$.  
\hfill\rule{1,5mm}{1,5mm}\\

\textbf{Lemma 5.2.} \textit{Let $G$ be an abelian $p$-group of rank 2 and order $p^n$, where $n\geq 3$. Then $\beta(G)\leq 1$ and the equality holds if and only if $G\cong\mathbb{Z}_2\times\mathbb{Z}_4$.}

\textbf{Proof.} Let $G\cong \mathbb{Z}_{p^{d_1}}\times\mathbb{Z}_{p^{d_2}}$ be an abelian $p$-group of order $p^n$, where $n\geq 3$, $1\leq d_1\leq d_2$ and $d_1+d_2=n$. There is only one abelian $p$-group of rank 2 and order $p^3$, this being $\mathbb{Z}_p\times\mathbb{Z}_{p^2}$. By (\ref{r6}), we have $\beta(\mathbb{Z}_p\times\mathbb{Z}_{p^2})=\frac{2p+4}{p^3}\leq 1$, and we note that the equality holds if and only if $p=2$. Let $n\geq 4$. By Lemma 3.5, it is sufficient to show that $\beta(\mathbb{Z}_{p^{\frac{n}{2}}}\times\mathbb{Z}_{p^{\frac{n}{2}}})<1$, if $n$ is even, and $\beta(\mathbb{Z}_{p^{\frac{n-1}{2}}}\times \mathbb{Z}_{p^{\frac{n+1}{2}}})< 1$, if $n$ is odd. We choose to analyse the second case. The same reasoning may be applied for the first one.  

Therefore, by (\ref{r6}), we obtain
\begin{align}\label{r8}
\beta(\mathbb{Z}_{p^{\frac{n-1}{2}}}\times \mathbb{Z}_{p^{\frac{n+1}{2}}})<1\Longleftrightarrow 2p^{\frac{n+3}{2}}-(n+3)p+n+1<(p-1)^2p^n.
\end{align}
But,
$$2p^{\frac{n+3}{2}}\leq (p-1)^2p^n\Longleftrightarrow 2\leq (p-1)^2p^{\frac{n-3}{2}},$$ which is true for   any odd integer $n\geq 4$ and any prime $p$. Consequently, inequality (\ref{r8}) holds and the proof is complete. 
\hfill\rule{1,5mm}{1,5mm}\\

\textbf{Lemma 5.3.} \textit{Let $p\geq 5$ be a prime number and let $G$ be an abelian $p$-group of rank 3 and order $p^n$, where $n\geq 3$. Then $\beta(G)<1.$}

\textbf{Proof.} Let $p\geq 5$ be a prime number. Let $G\cong\mathbb{Z}_{p^{d_1}}\times \mathbb{Z}_{p^{d_2}}\times\mathbb{Z}_{p^{d_3}}$ be an abelian $p$-group of order $p^n$, where $n\geq 3$, $1\leq d_1\leq d_2\leq d_3$ and $d_1+d_2+d_3=n.$ Using (\ref{r7}), one may prove that $\beta(\mathbb{Z}_p\times\mathbb{Z}_p\times\mathbb{Z}_{p^{n-2}})<1.$ Then, we may assume that $d_2\geq 2$ and $n\geq 5$.

According to (\ref{r7}), we must show that $p^{d_1+d_2+d_3}(p^2-1)^2(p-1)-A>0,$
where 
\begin{align*}
A= & (d_1+1)(d_3-d_2+1)p^{d_1+d_2+5}+2(d_1+1)p^{d_1+d_2+4}-2(d_1+1)(d_3-d_2)p^{d_1+d_2+3}-\\ &
2(d_1+1)p^{d_1+d_2+2}+(d_1+1)(d_3-d_2-1)p^{d_1+d_2+1}-(d_3+d_2-d_1+3)p^{2d_1+4}-2p^{2d_1+3}+ \\ &
(d_3+d_2-d_1-1)p^{2d_1+2}+(d_1+d_2+d_3+5)p^2+2p-(d_1+d_2+d_3+1).
\end{align*}
The main idea is to consider the left-hand side of the above inequality as a function of one variable $x$, which plays the role of $d_3$. In order to do this, we fix $n\geq 5$ and $d_2\geq 2$. Then $d_1=n-d_2-x$ and the constraint $1\leq d_1\leq d_2\leq d_3$ may be rewritten as $1\leq n-d_2-x\leq d_2\leq x$. So, we define the function $f:[d_2,n-3]\longrightarrow \mathbb{R}$ given by
\begin{align*}
f(x)=&[p^n(p^2-1)^2(p-1)+2(n-d_2-x+1)(x-d_2)p^{n-x+3}+2(n-d_2-x+1)p^{n-x+2}+n+1-\\&(n+5)p^2-2p-(n-d_2-x+1)(x-d_2-1)p^{n-x+1}-2(n-d_2-x+1)p^{n-x+4}-\\&(n-d_2-x+1)(x-d_2+1)p^{n-x+5}]+[(2x+2d_2-n+3)p^{2(n-d_2-x)+4}+2p^{2(n-d_2-x)+3}-\\&(2x+2d_2-n-1)p^{2(n-d_2-x)+2}].
\end{align*}
Our aim is to show that 
\begin{align}\label{r10}
f(x)>0, \ \forall \ x\in [d_2,n-3] \text{ such that } 1\leq n-d_2-x\leq d_2.
\end{align}

We remark that $f(x)=g(x)+h(x)$, where $g,h:[d_2,n-3]\longrightarrow\mathbb{R}$ are given by 
\begin{align*}
g(x)=&p^n(p^2-1)^2(p-1)+2(n-d_2-x+1)(x-d_2)p^{n-x+3}+2(n-d_2-x+1)p^{n-x+2}+n+1-\\&(n+5)p^2-2p-(n-d_2-x+1)(x-d_2-1)p^{n-x+1}-2(n-d_2-x+1)p^{n-x+4}-\\&(n-d_2-x+1)(x-d_2+1)p^{n-x+5},
\end{align*}   
and, respectively,
\begin{align*}
h(x)=(2x+2d_2-n+3)p^{2(n-d_2-x)+4}+2p^{2(n-d_2-x)+3}-(2x+2d_2-n-1)p^{2(n-d_2-x)+2}.
\end{align*}
Since $n-d_2-x\leq d_2$, we deduce that $2x+2d_2-n-1\geq x-1\geq d_2-1>0$. Then
$$(2x+2d_2-n+3)p^{2(n-d_2-x)+4}>(2x+2d_2-n-1)p^{2(n-d_2-x)+2}, \ \forall \ x\in [d_2,n-3],$$
and, consequently, it is obvious that $h$ takes only positive values. Also, we remark that $g$ is strictly  increasing since
\begin{align*}
g'(x)=& p^{n-x+1}(n-d_2-x+1)\lbrace -[(x-d_2)(2p^2-p^4-1)-(2p^4+2p^3-2p^2-2p)]\ln p\rbrace-\\& p^{n-x+1}[(x-d_2)(2p^2-p^4-1)-(p^4+2p^3-2p-1)]>0, \ \forall \ x\in [d_2,n-3].
\end{align*}
Then, to prove (\ref{r10}), it is sufficient to show that $g(d_2)\geq 0$, where
\begin{align*}
g(d_2)=&p^n(p^2-1)^2(p-1)+2(n-2d_2+1)p^{n-d_2+2}+(n-2d_2+1)p^{n-d_2+1}+n+1-\\& (n-2d_2+1)p^{n-d_2+5}-2(n-2d_2+1)p^{n-d_2+4}-(n+5)p^2-2p.        
\end{align*}
Now that we work with $x=d_2$, we note that the the constraint $1\leq n-d_2-x\leq d_2$ becomes $1\leq n-2d_2\leq d_2$. For our fixed variable $d_2\geq 2$, one can show that 
$$p^{d_2}>(1+d_2)(p-1), \ \forall \ p\geq 5.$$
Then, since $(p^2-1)^2(p-1)^2>p^5+2p^4, \ \forall \ p\geq 5,$ and $1+d_2\geq n-2d_2+1$, we have
$$p^{d_2}(p^2-1)^2(p-1)>(1+d_2)(p^2-1)^2(p-1)^2>(n-2d_2+1)(p^5+2p^4), \ \forall \ p\geq 5$$
and, consequently, we obtain that
\begin{align}\label{r11}
p^n(p^2-1)^2(p-1)>(n-2d_2+1)p^{n-d_2+5}+2(n-2d_2+1)p^{n-d_2+4}, \ \forall \ p\geq 5.
\end{align}
Going further, using Bernoulli's inequality, we get
$$2(n-2d_2+1)p^{n-d_2}\geq 2(n-2d_2+1)[1+(n-d_2)(p-1)]=2(n-2d_2+1)+2(n-2d_2+1)(n-d_2)(p-1).$$
But, $2(n-2d_2+1)>n-2d_2$ and 
$$2(n-2d_2+1)(n-d_2)(p-1)\geq 2(n-2d_2+1)(d_2+1)(p-1)>2(d_2+1)(p-1)>2d_2+5, \ \forall \ p\geq 5.$$
Then, $2(n-2d_2+1)p^{n-d_2}>n+5, \ \forall \ p\geq 5$, so
\begin{align}\label{r12}
2(n-2d_2+1)p^{n-d_2+2}>(n+5)p^2, \ \forall \ p\geq 5.
\end{align}
Finally, it is clear that 
\begin{align}\label{r13}
(n-2d_2+1)p^{n-d_2+1}>2p, \ \forall \ p\geq 5.
\end{align}
By adding the inequalities (\ref{r11}), (\ref{r12}) and (\ref{r13}), it follows that $g(d_2)>n+1>0$ and this completes our proof. 
\hfill\rule{1,5mm}{1,5mm}\\

Now, we may focus on the question that was marked at the beginning of this Section.\\

\textbf{Theorem 5.4.} \textit{Let $G$ be a finite abelian $p$-group. Then $\beta(H)\leq 1, \ \forall \ H\in L(G)$, if and only if $G$ is cyclic, $G$ is of rank 2 and $p\geq 3$, or $G$ is of rank 3 and $p\geq 5$.} 

\textbf{Proof.} Assume that $G$ is a finite abelian $p$-group such that $\beta(H)\leq 1, \ \forall \ H\in L(G)$. We want to show that $G$ has one of the above mentioned properties. Assume that $G$ has none of them. Then there are two possibilities: 
\begin{itemize}
\item[i)] $G$ has a subgroup $H$ such that $H\cong\mathbb{Z}_{2}^2$ or $H\cong\mathbb{Z}_{3}^3$;
\item[ii)] the rank of $G$ is at least 4. 
\end{itemize}
In the first case, we contradict our hypothesis since $\beta(\mathbb{Z}_2^2)=\frac{5}{4}>1$ and $\beta(\mathbb{Z}_3^3)=\frac{28}{27}>1$. In the second case, $G$ has a subgroup isomorphic to $\mathbb{Z}_{p}^4$ and, as we stated at the beginning of this Section, we have $\beta(\mathbb{Z}_{p}^4)>1$, a contradiction.

Conversely, if $G$ is cyclic, then all its subgroups are cyclic and, according to Lemma 5.1, we deduce that $\beta(H)\leq 1, \ \forall \ H\in L(G)$. If $G$ is of rank 2, $p\geq 3$ and $|G|=p^n$, where $n\geq 3$, then, by Lemmas 5.1 and 5.2, we conclude that $\beta(H)\leq 1, \ \forall \ H\in L(G)$. If $n=2$, then $G\cong\mathbb{Z}_{p}^2$ and the conclusion also holds since $p\geq 3$. In the end, if $G$ is of rank 3 and $p\geq 5$, then by our three preliminary results, it follows that $\beta(H)\leq 1, \ \forall \ H\in L(G),$ excepting the situation where $H\cong\mathbb{Z}_p^2.$ But, as we previously mentioned, the property holds for such subgroups. Consequently, the proof is complete.  
\hfill\rule{1,5mm}{1,5mm}\\

Using the multiplicativity of $\beta$, Theorem 5.4 may be extended to finite abelian groups.\\

\textbf{Corollary 5.5.} \textit{Let $k$ be a positive integer and $p_1,p_2,\ldots,p_k$ some distinct prime numbers. The only finite abelian groups $G$ satisfying $\beta(H)\leq 1, \ \forall \ H\in L(G),$ are $G\cong\xmare{i=1}k G_i$, where, for all $i\in\lbrace 1,2,\ldots,k\rbrace$, $G_i$ is a finite abelian $p_i$-group such that $G_i$ is cyclic, $G_i$ is of rank 2 and $p_i\geq 3$, or $G_i$ is of rank 3 and $p_i\geq 5$.}\\

We end this Section by marking a final result concerning the number of values that the function $\beta:L(G)\longrightarrow (0,\infty)$ may take, where $G$ is a finite abelian $p$-group of type $d=(d_k,d_{k-1},\ldots,d_1,0,\ldots)$, with $d_i\geq 1, \ \forall \ i \in\lbrace1,2,\ldots,k\rbrace$. A consequence of the proof of Lemma 5.1 is the fact that $\beta$ is injective on  $L_1(G)$. Then, it follows that $|Im \ \beta|\geq d_i , \ \forall \ i\in\lbrace 1,2,\ldots,k\rbrace$, and $|Im \ \beta|$ is greater than or equal to the number of different positive integers that form the partition $d$.\\

Using the above ideas, one can prove the following result.\\

\textbf{Theorem 5.6.} \textit{Let $G$ be a finite abelian $p$-group. Then:
\begin{itemize}
\item[i)] $|Im \ \beta|=1\Longleftrightarrow G\cong\mathbb{Z}_2.$
\item[ii)] $|Im \ \beta|=2 \Longleftrightarrow G\cong \mathbb{Z}_p$, where p is odd, $G\cong\mathbb{Z}_3^2, G\cong\mathbb{Z}_4$ or $G\cong\mathbb{Z}_2^2.$ 
\end{itemize}}      
   
\section{Further research}

In Section 3, we saw that $\beta$ is not monotonic on the class of finite abelian $p$-groups of a given order. Also, in Section 2, we provided an example to show that $\beta$ is not injective on the class of finite $p$-groups of a given order.\\ 

\textbf{Problem 6.1.} \textit{Is $\beta$ injective on the class of finite abelian $p$-groups of a given order? In particular, suppose that two abelian $p$-groups $G_1$ and $G_2$ or order $p^n$ have the following chains of subgroups
$$\lbrace1\rbrace=H_0\leq H_1\leq\ldots\leq H_{n-1}\leq H_n=G_1 \text{ and } \lbrace 1\rbrace=K_0\leq K_1\leq\ldots\leq K_n=G_2,$$
where $|H_i|=|K_i|=p^i$ and $\beta(H_i)=\beta(K_i), \ \forall \ i\in\lbrace 0,1,\ldots,n\rbrace$. It is true that $G_1\cong G_2?$}\\

There are some classifications of finite groups $G$ that have $|G|-r$ cyclic subgroups, where $0\leq r\leq 5$. For more details, the reader may refer to \cite{2,25}. An interesting problem would be to classify all finite groups $G$ satisfying $\beta(G)=1$. This is equivalent to finding all finite groups $G$ for which $|G|=|L(G)|$. A starting point would be to work only with finite abelian $p$-groups. In this case, some additional comments can be written. Hence, let $G$ be a finite abelian $p$-group of rank $k$ and order $p^n$. In Section 5, we proved that if $k\in\lbrace 1,2,3,4\rbrace$, then $\beta(G)=1$ if and only if $G\cong\mathbb{Z}_2$ or $G\cong\mathbb{Z}_{2}\times\mathbb{Z}_4$. Hence, we may assume that $k\geq 5$. Obviously, $G$ has a subgroup $H\cong\mathbb{Z}_p^k$. By Theorem 2.1 of \cite{1}, we have $p^{\frac{k^2}{4}}<|L(H)|<p^{\frac{k^2}{4}+4}$. Then, if $k>2\sqrt{n}$, it follows that $\beta(G)>1.$ Hence, to find additional abelian $p$-groups $G$ of rank $k$ and order $p^n$ having the property $\beta(G)=1$, one may assume that $5\leq k\leq 2\sqrt{n}$.     
\\

\textbf{Problem 6.2.} \textit{Classify all finite abelian $p$-groups $G$ that satisfy $\beta(G)=1.$}\\

$\beta$ viewed as a function on the subgroup lattice of a finite abelian $p$-group is not monotonic, in general. For instance, if $G\cong\mathbb{Z}_2\times\mathbb{Z}_4$, by taking two subgroups $H_1\cong\mathbb{Z}_2$ and $H_2\cong\mathbb{Z}_2^2$, we have $H_1\subset H_2\subset G$, but $\beta(H_1)=\beta(G)=1$, while $\beta(H_2)=\frac{5}{4}$. Hence, the following question arises.\\

\textbf{Problem 6.3.} \textit{Classify all finite abelian $p$-groups $G$ for which $\beta:L(G)\longrightarrow(0,\infty)$ is monotonic.}\\

For a finite group $G$, the quantity $\alpha(G)=\frac{|L_1(G)|}{|G|}$ was introduced in \cite{13}, the authors obtaining a lot of results by studying this ratio corresponding to different classes of finite groups. Obviously,  
$$\alpha(G)\leq \beta(G),$$
and the equality holds if and only if $G$ is cyclic.\\

\textbf{Problem 6.4.} Study the connections between $\alpha$ and $\beta$ on a class of finite groups.

\vspace*{3ex}
\small
\hfill
\begin{minipage}[t]{6cm}
Mihai-Silviu Lazorec \\
Faculty of  Mathematics \\
"Al.I. Cuza" University \\
Ia\c si, Romania \\
e-mail: {\tt Mihai.Lazorec@student.uaic.ro}
\end{minipage}

\begin{thebibliography}{100}
\bibitem{1} Aivazidis, S., {\it The subgroup permutability degree of projective special linear groups over fields of even characteristic}, J. Group Theory {\bf 16} (2013), 383-396.

\bibitem{2} Belshoff, R., Dillstrom, J. and Reid, L., \textit{Finite groups with a prescribed number of cyclic subgroups}, to appear in Comm. Alg., https://www.tandfonline.com/doi/full/10.1080/00927872.2018.1499923.

\bibitem{3} Berkovich, Y., \textit{Groups of prime power order}, vol. 1, de Gruyter, Berlin, 2008.

\bibitem{4} Bhowmik, G., \textit{Evaluation of divisor functions of matrices}, Acta Arith. \textbf{74} (1996), 155-159.

\bibitem{5} Birkhoff, G., \textit{Subgroups of Abelian groups}, Proc. London Math. Soc. \textbf{38} (1934-35), 385-401.

\bibitem{6} Borovik, A. V., Pyber, L. and Shalev, A., \textit{Maximal subgroups in finite and profinite groups}, Trans. Amer. Math. Soc. \textbf{348} (1996), 3745-3761.

\bibitem{7} Butler, L. M., \textit{A unimodality result in the enumeration of subgroups of a finite abelian group}, Proc. Am. Math. Soc. {\bf 101} (1987), 771-775.

\bibitem{8} Butler, L. M., \textit{Subgroup Lattices and Symmetric Functions}, Mem. Amer. Math. Soc. {\bf 112}, 1994.

\bibitem{9} C\u alug\u areanu, G., \textit{The total number of subgroups of a finite abelian
group}, Sci. Math. Jpn. \textbf{60} (2004), 157-167.

\bibitem{10} Chen, Y. and Chen, G., \textit{Finite groups with the set of the number of subgroups of possible order containing exactly two elements}, Proc.: Math. Sci. \textbf{123} (2013), 491-498.

\bibitem{11} Delsarte, S., \textit{Fonctions de M\" obius sur les groupes abeliens finis}, Annals of Math. {\bf 49} (1948), 600-609.

\bibitem{12} Frobenius, G., \textit{Verallgeminerung des Sylow'schen Satzes}, Sitzungsberichte
der K\" oniglich Preu\ss ischen Akademie der Wissenschaften zu Berlin, 1895, 981-993.

\bibitem{13} Garonzi, M. and Lima, I., \textit{On the number of cyclic subgroups of a finite group}, Bull. Brazil. Math. Soc. \textbf{49} (2018), 515-530.

\bibitem{14} Hampejs, M. and T\'{o}th, L., \textit{On the subgroups of finite abelian groups of rank three}, Annales Univ. Sci. Budapest., Sect. Comp. \textbf{39} (2013), 111-124.

\bibitem{15} Kulakoff, A., \textit{\" Uber die Anzahl der eigentlichen Untergruppen und der Elemente von gegebener Ordnung in p-Gruppen}, Math. Ann. \textbf{104} (1931), 778-793.

\bibitem{16} Macdonald, I. G., \textit{Symmetric Functions and Hall Polynomials}, 2nd ed. Oxford Science Publications. Oxford: Clarendon Press, 1998.

\bibitem{17} Nitecki, Z., \textit{Cantorvals and subsum sets of null sequences}, Amer. Math. Monthly \textbf{122} (2015), 862-870.

\bibitem{18} Oh, J. M., \textit{An explicit formula for the number of subgroups of a finite abelian p-group up to rank 3}, Commun. Korean Math. Soc. \textbf{28} (2013), 649-667.

\bibitem{19} Qu, H., Sun, Y. and Zhang, Q., \textit{Finite p-groups in which the number of subgroups of possible order is less than or equal to $p^3$}, Chin. Ann. Math. Ser. B \textbf{31} (2010), 497-506.

\bibitem{20} Qu, H., \textit{Finite non-elementary abelian p-groups whose number of subgroups is maximal}, Israel J. Math. \textbf{195} (2013), 773-781.

\bibitem{21} R. Schmidt, {\it Subgroup lattices of groups}, de Gruyter Expositions in Ma\-the\-ma\-tics 14, de Gruyter, Berlin, 1994.

\bibitem{22} M. Suzuki, {\it Group theory}, I, II, Springer Verlag, Berlin, 1982, 1986.

\bibitem{23} M. T\u arn\u auceanu, {\it Subgroup commutativity degrees of finite groups}, J. Algebra {\bf 321} (2009), 2508-2520.

\bibitem{24} T\u arn\u auceanu, M., {\it An arithmetic method of counting the subgroups of a finite abelian group}, Bull. Math. Soc. Sci. Math. Roumanie (N.S.) {\bf 53} (2010), 373-386.

\bibitem{25} T\u arn\u auceanu, M., \textit{Finite groups with a certain number of cyclic subgroups},
Amer. Math. Monthly \textbf{122} (2015), 275-276.

\bibitem{26} T\u arn\u auceanu, M. and T\'{o}th, L., \textit{On the number of subgroups of a given exponent in a finite abelian group}, Publ. Inst. Math. (Beograd) (N.S.) \textbf{101} (2017), 121-133.

\bibitem{27} T\u arn\u auceanu, M., \textit{On a conjecture by Haipeng Qu}, accepted for publication in J. Group Theory, https://arxiv.org/abs/1811.07478.

\bibitem{28} The GAP Group, \textit{GAP -- Groups, Algorithms, and Programming, Version 4.9.3}, 2018, https://www.gap-system.org.

\bibitem{29} Yeh, Y., \textit{On prime power abelian groups}, Bull. Amer. Math. Soc. {\bf 54} (1948), 323-327.
\end{thebibliography}
\end{document}